\documentclass[12pt, letterpaper]{article}

\usepackage{wasysym}
\usepackage{booktabs}
\usepackage[T1]{fontenc}
\usepackage{graphicx}
\usepackage[latin,english]{babel}
\usepackage[noeledmac]{ledmac}
\usepackage{ledpar}
\lineation{section}
\usepackage{graphicx,vmargin}
\usepackage{url}
\setpapersize{USletter}

\setmarginsrb{0.8in}{0.8in}{0.8in}{0.9in}{0pt}{0pt}{0pt}{24pt}

\begin{document}

\begin{center}
\textbf{\begin{Large}\textit{Mensurae Universales Magnitudinum ac Temporum} \\by Adam Adamandy Kocha\'nski -- Latin text with annotated English translation      
\end{Large}}
\end{center}
\begin{center}
translated by Henryk Fuk\'s\\
Department of Mathematics and Statistics,
Brock University, St. Catharines, ON, Canada\\
email \texttt{hfuks@brocku.ca}\end{center}
\vspace{0.5cm}
\noindent \textit{Translator's note:} The Latin text of \textit{Mensurae Universales} presented here closely  follows the original text published in 
\cite{mensuniv1} (republished in \cite{1740opuscula}). Punctuation, capitalization, and mathematical  notation have been preserved. 
The translation is as faithful as possible, often literal, and it is mainly intended to be of help to  those  who wish to study the  
original  Latin text. Footnotes and figures in the appendix have
been added by the translator.

\vspace{0.3cm}

\setlength{\Lcolwidth}{0.465\textwidth}
\setlength{\Rcolwidth}{0.495\textwidth}

\numberlinefalse

\begin{pairs}
\begin{Leftside}
\beginnumbering
\selectlanguage{latin}
\noindent\setline{1}
\pstart
\begin{center}
\emph{\begin{large}ADAMI ADAMANDI KOCHANSKI\\E SOCIET. JESU  \end{large}\\
Sereniss. Poloniarum Regis Mathematici,\\
MENSURAE Universales Magnitudinum ac Temporum\\
Ad Actorum Erud. Collectores.
}\end{center}
\mbox{\,}
\pend

%
%

\pstart
Praeterlapsi anni 1685 Mense Decembri, transmiseram 
Parisios meas aliquot Observationes  ad Amicum quendam,
rogatum impense, ut eas ipsemet Ephemeridum Gallicarum
Collectori in manus traderet: sed eum trimestri iam elapso
nihil, comperti habeam, utrum eo pervenerint, visum mihi est
e re literaria futurum, si eas, prout a me nunc conscriptae
fuerant, ab interitu vindicarem. Omissa itaque conjectura illa
mea de Metu Diurno Telluris, e terra motu petita, quae 
Actis Erud. Lipsiensibus Mense Julio Anni 1685 inserta est, \&
quam in Gallia, succussionibus ejusmodi magis, quam nostrae
regiones obnoxia, notam reddere, \& observatores exstimulare
cupiebam, reliqua eo tunc transmissa, verbis iisdem, fideliter
\& integre, exscribam, ut sequitur.
\pend
\pstart
\emph{Claris. \& Eruditis. Viro Ephemeridum Gallic. Collectori \& c.}

Non aliud paucis his, ut Tibi nunc, Cl. Vir, adfero,
praetexam patrocinium praefationemque, quam nos atque natura
sciendi cupidos, ac ad ea, quas scimus aliis enarranda proclives
esse. Eo quoque invitant Tua, aliorumque in promovendis
scientiis exempla. Quod si non omnes ad exstruenda
scientiarum palatia, marmor aut gemmas adferimus, etiam 
sabulum, paleaque leves caemento hujus aedificii non erunt, 
opinor, inutiles.
\pend

\pstart
\begin{center}
I \emph{Nova mensura universalis Temporum}. 
 \end{center}

Quaerit hucusque Geometria certam quandam, ac invariabilem
Corporum Linearumque mensuram; quae quidem, opinione
mea, in ipsius magnitudinis continuae natura delitescit, sed
non prius eruenda, quam cognita sit, \& aliquam mensuram
revocata quantitas Angulorum Mixtilineorum, quales sunt
semicirculorum, segmentorumque: At quia id vicissim certam
lineae rectae, arcus eos subtendentis, quantitatem praesupponit,
nec fortassis alia via indagari potest, difficultas  haec abit in 
circulum nullo humano ingenio explicabilem. Quocirca Problemati
huic solvendo adhibita fuit a nonnullis mensura Temporis a
Pendulo petita, cuius vibrationes minuta secunda primi Mobilis,
vel horum partes aliquotas adaequant. Et vero improbari
minime debet Inventum istud, si debitis cum cautelis in eo
procedatur, ac praesertim ut nobis certo constet, Pendula nec
ab Athmosphaerae in diversis climatibus varietate, nec a varia
centri terrae distantia sensibiliter alterari; quod an quisquam
exploraverit, ignoro: quia tamen ejusmodi penduli inventio
hominibus  occupatis admodum operosa videtur: nam plurium
vicariam opem postulat, unico tentamine non obtinetur, nubibus
eluditur \&c. Idcirco innuam aliud fortasse commodius,
cujus ope motus omnes etiam Caelestes mensurari, ac praedictorum
Pendulorum vera quantitas facilius determinari poterit.
Observavi ab annis compluribus, candelas usitate cujusvis
crassitiei, ardentes in aere libero ac tranquillo, scintillationes,
sive subsultationes illas suas ab aere undulatim accurrente 
ortas, vibrationibus ad sensum aequidiuturnis peragere. Inveniatur
ergo Pendulum, cuius una vibratio quaternis, aut fenis
ejusmodi palpitationibus sit aequalis, atque hujus Penduli 
longitudo simpla vel dupla cum motu diurno fixarum, nec non
solis medio, conferatur: sic enim habebimus mensuram non e
Caelo, sed ex Elemento ignis nostri sublunaris deductam;
quam tamen \& ab aliis elementis aliter obtinebimus.
\pend
\pstart
II. \emph{Nova mensura universalis Magnitudinum.}

Praeter aliam indicatam, est alia quadam non inelegans
huius indagandae ratio, ab elemento Aquae deducta, diversa 
tamen ab ea, quam T. Livius Burattini, amicus quondam meus,
secutus est in suo hac de re tractatu, Italice in Lituania 
edito: sed quia meam hanc Observationem paucis hic explicare
non possum, alteri tempori locoque reservare cogor. Nostrae
igitur Arti Geometriae Practicae, a naturis Aereis mensuras 
mutuabimur hoc modo. Praetermissis itaque omnium vegetabilium
seminibus, pilis animalium, fibris corporum Microscopio 
detectis \&c. ego spem omnem collocavi in Volucribus. Aptissimas
huic negotio crederem pennas ex alis passeris domestici ruralis,
supposito tamen quod haec species avicularum etiam in extremo
Oriente, nec non America reperiatur, \& quidem ejusdem staturae
cum nostratibus. Refert O. Dappers in sua Africa Germanicae
edita, passeres ad caput Bonae spei reperiri; sed earum pennae
cum nostratibus essent conferendae. Praxis ad rem praesentem ita
esset instituenda. Ex alterutra ala passeris masculi \& senis, sumatur
certa quaedam penna, v. g. tertia \& quarta, \&  ea primum 
artificio facili in directum extendatur: tum in medio cauliculi, 
plumulis transversis convestiti, capiatur certus quidam numerus 
intervallorum, quibus transversae plumulae ab invicem distant, istae
enim divisiones inservient nobis per modum scalae cujusdam artificialis,
juxta nam taxari poterant partes aliquotae pollicum, aut
palmorum pedis Geometrici, ac etiam aliarum Urbium, ac posteritati
scripto commendari: non enim mihi credibile videtur, passerculorum
nostrorum proceritatem secuturis saeculis immutandam.
Cautelas in operando necessarias non enumero, quia perspicacibus
loquor. Si quis autem subtiliore, \& exactiore scala uti volet, hanc
inveniet in iisdem passerum pennis: nam \&  ipsaemet transversae
plumulae subdivisae sunt aliis transversis plumulis, quae iudice
rudiore Microscopio, magis aequalibus ab invicem distant intervallis.
\emph{Inquirant igitur curiosi}, quibus aut maria trajicere, aut a regionum
peritis id intelligere commodum fuerit , utrum passeres
\& hirundines, aut alia aves communes ubique terrarum , \&  quidem
( quod caput rei est) ejusdem cum nostris staturae reperiantur.
Plumulas transversas struthionis, pavonis, ciconiae \&c. 
relinquendas censeo, quod non sint adeo obviae, quales esse praestat.
Habeo complures alas passerum e diversis Europae partibus,
qua primo adspectu aequales esse videntur, sed eas per otium diligentius
adinvicem conferam. 
\emph{Pondus universale} tunc primum adferam, 
cum explicuero artificium illud aquam tractandi initio commemoratum.
Id interim silentio minime praeteribo, pro pondere
universali rudioris minervae, sumi posse aureos Hungaricos, quos
\&  Ducatos appellamus: hi enim per totum orbem sparsi sunt , \& 
quidem ponderis delectu prius instituto, non admodum inaequalis.
His igitur utcunque taxari poterunt librae (pondo) diversarum urbium,
etiam in Asia, ubi Anglica Marcarum ponduscula ad manum haberi non possunt.
\pend
\pstart
III. \emph{Mensurarum facilis ad posteros, vel absentes transmissio.}

Non abs re fortassis erit hoc loco commemorare, deceptum fuisse
W. Snellium, \& aliquos cum eo Auctores, qui crediderunt mensuras
pedum chartae madidae impressas, certa quadam proportione
post ejus siccationem contrahi. Certum enim est, chartae longitudinem
magis huic contracturae subesse, quam latitudinem.
\emph{Longitudinem}
autem illam voco, iuxta quam porriguntur in longum 
vestigia filorum aeneorum, e quibus cribrum aut forma papyri contexta
est. \emph{Latitudo chartae ea est}, quae praedicta filamenta transversim
connectit geminis filis, hinc  inde ad intervallum pollicis circiter.
Occasione porro pedis Romani antiqui, a Villalpando in Apparatu
Templi chartae impressi, \& genuino pedi plerumque congruentis,
exploravi diversas in Europa chartas, \& inter eas, illam ipsam Hispanicam Villalpandi, quae inter omnes minimam passa est a siccatione
contracturam secundum latitudinem, in quam Pes ille porrigitur:
fortassis autem id accidit ex eo, quod ipsa a quidem tenuis fit, licet 
cribro crassioribus filis contexto formata fuerit; nam in alia quadam 
ejusdem fortis, at nonnihil crassiore, inveni totum Pedem una centesima
breviorem eo, quem Romae a Grunbergero nostro, Villalpandi,
in calculandis Apparatus Tabulis, adjutore, lapidi subtiliter incisum,
exactissime dimensus sum, \&  in aenea lamina designatum adhuc 
possideo. Hic tamen alter Pes impressus crassiori chartae, nec
dum compactoris malleum subierat, quem postea passus, defectum
illum liberaliter compensavit.

Experimentum illius rei sic institui: e qualibet chartarum infra
positarum, sumpti duas fascias latas 130 centesimis unius unciae pedis
Rhenani, longas autem Semipede Rhen. quarum fasciarum una erat
resecta juxta longitudinem, altera secundum latitudinem chartae.
Hae fasciae aqua fluviatili madefactae creverunt juxta numeros
adjectos, qui sunt centesimae partes unius unciae Pedis Rhinlandici.
\pend
 
\end{Leftside}

\begin{Rightside}
\beginnumbering
\selectlanguage{english}
\noindent\setline{1}
\pstart
\begin{center}
\emph{\begin{large}ADAM ADAMANDY KOCHA\'NSKI \\FROM THE SOCIETY OF JESUS, \end{large}\\
Mathematician of the most serene King\footnote{John III Sobieski (1629 -- 1696), from 1674 until his death King of Poland and Grand Duke of Lithuania.} of Poland,\\
Universal MEASURES of Magnitude and Time\\
To the Editors of Acta Eruditorum
}\end{center}
\mbox{\,}
\pend

\pstart
In the month of  December of the passing year 1686
I forwarded to Paris to a certain friend   a few of my observations, asking him zealously to deliver them in person to the hands of the editor of Ephemerides Eruditorum Gallicae\footnote{\emph{Journal des Sçavans}, academic journal estab. in 1665 by Denis de Sallo, French writer and lawyer.}; yet since a trimester already passed,
and I do not know for sure whether they arrived,
it seems to be proper, for the sake of the future of
the scholarly cause, to save them from perishing.
Omitting, therefore, my conjecture regarding the
daily motion of the Earth, derived from [the  observations of] earthquakes, which was inserted
in Acta Eruditorum in the month of July of the year 1685, and which I desired to be noted and to stimulate observers in France, liable to stronger shakes than our regions, I will copy the rest of what was then submitted, with the same wording,
faithfully and wholly, as follows.
\pend
\pstart
\emph{To the most illustrious and learned 
man\footnote{Fr. Jean-Paul de La Roque (?--1691), editor of
\emph{Journal des Sçavans} in the years 1674-1687.}, editor of Ephemerides Eruditorum Gallicae etc.}

I will cover the littleness of that which I send you with no other defense and
introduction than this: by nature we desire knowledge, and we are inclined
to pass what we know on to others. Your examples of the advancement of science, as well as the examples of others, also encourage this. Although not everybody contributes marble and precious stones  to the construction of the palace of sciences, even
gravel and light straw are not, I suppose, inexpedient as components of
the concrete [used in construction] of this building. 
\pend

\pstart
\begin{center}
I. \emph{New measure of universal Time}. 
\end{center}

Geometry searches thus far for a fixed and invariant measure of bodies and lines;
this, however, in my opinion, nature conceals in its own magnitudes,
and it is to be extracted not sooner then it is related and linked to a
measure of quantity of mixtilinear angles, which are of the kind of semicircles and segments. And because this presupposes a certain quantity of the straight line, 
stretching beneath curves, and because perhaps another way cannot be found, this difficulty leads to a circular problem, not resolvable by human ingenuity.

On the account of solving this problem, many employed a measure of time
derived from a pendulum, whose oscillations are equal to seconds of time
of the prime motion, or some parts thereof. Indeed, this invention should not be
rejected if one proceeds with due precautions, and particularly if
we make sure that the pendula are not perceptibly affected by differences of the atmosphere in different climates, nor by differences in the distance from the center of the earth.
I do not know if somebody explored this issue. 
Still, this method seems to be troublesome for the people using it: it requires the effort of several helpers,
the result is not obtainable by a single attempt, it is frustrated by clouds, etc.
On that account, I will hint another, perhaps more convenient method, by the power of which the motion of all,
even celestial bodies, can be measured, and for which the true quantity can be determined easier than with 
the aforementioned pendulum. I have observed for many years commonly used candles  of any whatever thickness,
burning in the  unencumbered and calm air, their twinkling or flickering arising from the air rushing in
a wavelike fashion, causing oscillations perceived as of equal duration. Let a pendulum be contrived,
such that its period is equal to four or six such oscillations, and let its single or double length be compared to  
the daily motion of fixed [stars], or of the middle of the sun. In such a way we will have a measure derived not
from the sky, but from the element of our sublunar fire; we will yet obtain it differently  also from other elements.
\pend
\pstart
II. \emph{New universal measure of magnitude}

Before I point out another, there is also a certain other [measure] not inelegant by its investigation method,
derived from the element of water, nevertheless separate from that which was pursued by formerly a friend of my 
T. Livius Boratini\footnote{Tito Livio Burattini (1617--1681),
Italian nobleman, inventor, architect and diplomat, from 1641 till his death working in the service of Polish kings.
Author of
\emph{Universal measure} (orig. \emph{Misura universale}, Vilnius 1675).}, in the treaty on this subject, published in Lithuania in Italian; but because I cannot explain 
this observation of mine in a few words, I am compelled to reserve it for another time and place.
From the nature of air, therefore, we will move the measures to the practice of our art of geometry.
Passing over the seeds of all plants, animal hairs, fibers of materials revealed by a microscope etc.,
I placed all hope in flying creatures. Most suitable for this work I woul believe to be feathers from the 
wings of domestic sparrow from the countryside.
Supposedly this species of bird is found even in the far East in America, and indeed its stature is the same 
as the stature of our birds. O. Dappers\footnote{Olfert Dapper (1636 -- 1689), Dutch physician and writer, author of \emph{Description of Africa} (orig. \emph{Naukeurige Beschrijvinge der Afrikaensche Gewesten}, Amsterdam 1668).}
reports in his ``Africa'', published in German, that sparrows are found at the 
Cape of Good Hope, and its feathers are comparable [with the feathers of our sparrows].
Practice pertaining to our matter is established as follows. From one of the two wings of a masculine and aged sparrow
a certain feather is taken, for example the third and the fourth one, and at first  it is skilfully made straight.
Then in the middle of the stem vested with transverse barbs, we take a  number of intervals by which
transverse barbs are separated from each other, and this will serve  us as a unit of an artificial scale
by which one can express the length of fractions of inches or palms entering geometric foot, and feet from
others cities, and also entrust it to writing for posterity: it does not seem credible to me that
the size of the sparrow would change in the coming ages.
I do not specify precautions necessary in this operation, because I speak to acute [readers].
If somebody,  moreover, wanted to use a more subtle and more exact scale, he will find it in the 
feathers of the same sparrow, for the transverse barbs are subdivided by other transverse barbules which,
judging by [observations made with] a crude microscope, stand apart from each other with even more equal intervals.
Let, therefore, \emph{the inquisitive}, for whom it is opportune to cross the sea or to 
gain this knowledge from local experts, \emph{examine} whether sparrows and swallows or other birds common 
to all lands indeed (which is crucial to our considerations) share the same stature with our birds.
The transverse barbs of ostriches, peacocks, storks, etc. I consider to be better left behind, because 
it is not obvious so far how good they are [for our purpose].
I have many wings of sparrows form different parts of Europe, which at first sight appear to be equal,
but I will compare them more diligently in my free time.
Let me report on the \emph{universal weight} at the time when I will be explaining the method
of water pumping mentioned at the beginning. In the meanwhile I will not pass over this completely silently:
for a crude universal weight one can take Hungarian gold pieces which we also call ducats.
These indeed are scattered all over the world and certainly not exceedingly different from 
previously instituted choices of weights.
With them, therefore, one can asses the weight of the pound used in different cities, even in Asia,
where English Mark\footnote{Marca Anglicana, weight of one English pound used with balance scales.} weights  cannot be obtained. 
 \pend
 \pstart
 III. \emph{Ease of Transmission of measures to posterity or those absent}
 
 Perhaps it is not without merit to mention in this place 
 that W. Snellius\footnote{Willebrord Snellius (1580-1626), Dutch astronomer and mathematician.}, and with him other authors, have been deceived
 believing that measures of foot printed on wet sheets of paper, contract
 proportionally after 
 drying\footnote{In his book \emph{Eratosthenes Batavus} (Logundi Batavorum, 1617), Snellius states that
 \emph{forma chartae impreffa sexagesimam partem ab archetypo suo deducat}, ``a shape printed on paper takes away 1/60 
 from its pattern''. He assumes it is the same in all directions.}.
 It is, in fact, certain that paper
 will contract more in longitude  than in latitude. By \emph{longitude}
 I mean the direction of traces of bronze filaments from which
 the sieve or form of the paper is made. The  \emph{latitude} is the direction
 of transverse twin filaments connecting the aforementioned ones, going back and forth 
 at intervals of about the width of the thumb.
 Hereafter under the pretext of investigating the ancient Roman foot,
 printed on the page\footnote{See Fig. \ref{vasefig} in the Appendix.} of the book ``Equipment of the Temple'' by Villapandus\footnote{Juan Bautista Villalpando SJ (1552 --1608), Jesuit scholar, mathematician, and architect,
 author of \emph{Ezechielem Explanationes} in which he reconstructed the Temple of Solomon. The third volume of this work was published as 
\emph{Apparatus Urbis ac Templi Hierosomitani}, Romae 1604.},
 and mostly congruent to the genuine Roman foot, I explored various European papers,
 among them the the same Spanish paper of Villapandus, on which this foot
 is printed: among others it suffers the least contraction in the latitude while
 drying. Perhaps it happened because it was thinner than the others, and yet
 formed with sieve weaved with more coarse filaments.
 But in others equally strong and a bit thicker I found the foot to be
 one hundredth shorter from the subtly cut stone pattern of the Roman foot
 made in Rome by our Grunberger\footnote{Christoph Grienberger SJ (1561 -- 1636), Austrian 
Jesuit mathematician and astronomer, author of a catalog of fixed stars.
}, helper of Villpandus in the calculations
 of the tables for ``Apparatus''; I have measured this pattern most precisely and now I have 
 it marked on a brass plate. Nevertheless, another foot printed on thicker paper,
 not subjected to the mallet of the press man, which it suffered afterwards,
 compensated this defect graciously.

 I instituted this experiment as follows: from whichever sheet of paper
 shown below, I took two streaks 130/100 of 1/12 of the Rhineland foot\footnote{Rheinfuss, unit of length equal to 31.4cm.}
 wide,
 1/2 foot long, of which one was    cut along the longitude, the other along the
 latitude of the sheet. These streaks, made wet with river water,
 grew according to the numbers added [to the table], which are hundredths of
 1/12 of the Rhineland foot.
 \pend

\end{Rightside}

\Columns

\end{pairs}

\vskip 1cm

\begin{tabular}{lr|r}
\emph{Sequentes Chartae madefactae creverunt} &\emph{In Longum} &\emph{In Latum} \\
Romana di Fuligno candida \& solida     &10  &7.\\
Hispanica Villalpandi in figura Congii Rom.   &10  &5.\\
Eadem nonnihil  crassior, de qua superius sermo  &10  &7.\\
Florentina epistolaris, quales \& septem consequentes &13  &8.\\
Genuensis cum tridentis  insigni      &16 &9.\\
Italica quaedam tenuissima, cum hasta    &17  &10.\\
Gallica cum Cornu Tabellionis      &10  &9.\\
Hollandica praecedentem mentita     &11 &8. \\
Joachimicae Vallis in Bohemia, candida, tenuis, at
 compacta         &20 &13.\\
Commotoviensis in Bohem. probe aluminata, tenuis 
tamen         &14  &10.\\
Uratislaviens. in Silesia, tenuis, at compacta  &13  &10.\\
Regalis crassa Bohemica      &14  &9. \\
Dantiscana crassa cum Carpionis insigni    &9  &7.\\
ltalica Regalis vetusta, ac solide compacta   &10  &9.\\
Alia crassis filamentis contexta, \& mediocriter crassa &13   &9 1/2 
\end{tabular}

 \begin{tabular}{l|r|r}
\emph{Following wet sheets have grown} &\emph{in lenght} &\emph{in width} \\
Roman from Foligno bright  \& firm    &10  &7.\\
Spanish of Villalpandi with the figure of Roman vase   &10  &5.\\
The same a bit  thicker, of which is spoken above  &10  &7.\\
Letter paper from Florence, likewise also seven following &13  &8.\\
From Genua with the sign of trident       &16 &9.\\
Italian most thin, with spear    &17  &10.\\
French with messenger's horn      &10  &9.\\
Dutch imitating the previous one     &11 &8. \\
Of Joachim Vallis in Bohemia, bright, thin, and
 close-packed        &20 &13.\\
Of Chomutov  in Bohemia well treate with alun, yet thin          &14  &10.\\
Of Wroc{\l}aw/Breslau in Silesia, thin and close-packed  &13  &10.\\
Royal Bohemian thick      &14  &9. \\
From Gda{\'n}sk/Danzig thick with sign of Carpio(?)    &9  &7.\\
Italian Royal old,  firmly compacted   &10  &9.\\
Another weaved with thick filaments, \& moderately thick &13   &9 1/2 
\end{tabular}

\begin{pairs}
\begin{Leftside}
\selectlanguage{latin}
\noindent\setline{1}
\pstart
Observavi praeterea, lineas ac figuras laminis aeneis incisas,
\& usitato praelo Cylindrico in charta expressas, magis ea parte
contrahi, qua Cylindrorum violenta pressione charta fuerat 
in longum distenta, quod non ita contingit in transversum;
unde fit, ut circuli praesertim majusculi, laminae diligenter 
incisi, \&  memorato prelo expressi, charta siccata reddantur 
altera parte contractiores. Sed neque Ricciolus noster in sua
Geogr. Reform. errandi periculum omne sustulit, semipedem
Romanum chartae siccatae imprimi curando; nam haec ipsa ad
aeris alterationes nonnihil mutatur, \&  tractu temporis non 
secus ac tabularum lignearum, quin \&  eboris latitudo, minuitur.
Securior itaque modus esset hujusmodi mensuras particulares ad
posteros transmittendi, incidendo illas formae metallicae, ac tum
imprimendo bracteolis illis aeris coronarii, quae ab Italis 
\emph{Orpello}, a Germanis \emph{Rauschgold} appellantur. Hae enim non 
incommode libris a compactore inferi poterunt. Alteratio illa 
modica, quam subeunt a calore vel frigore, sensibilem in semipedes
non inducet errorem. Si quis tamen e soliis impressis genuinas
consequi volet mensuras, eas a charta prius humectata, \&  
convenienter extensa, sat exactas poterit impetrare.

  Per Epistolam vero commodissime mitti poterit Amico 
absenti, mensura semipedis, vel ejus quadrantis, in chorda tenui
aenea, quales in Musicis instrumentis adhibentur; In ejus enim
extremis ad intra reflexis fiunt per contorsionem chordae duo
laquei A \&  B, partes vero contortae malleo nonnihil tunduntur, 
ut firmius cohaereant: tandem extento supra datam 
mensuram mediocritier filo aeneo, eius extremitates internae 
forficula, vel cultello, accurata manu quantum opus est, reciduntur.
Inspice adjectam Figuram.  
\pend
 
\end{Leftside}

\begin{Rightside}
\selectlanguage{english}
\noindent\setline{1}
 \pstart
I observed thereafter that lines and figures engraved in copper plates, when they are printed on paper using cylindrical press,
contract more in the direction in which the sheet was stretched along by the violent pressure of the cylinders; this does not happen so
in the transversal direction. Because of this, it happens that circles, particularly larger ones, diligently engraved in the plate
and printed with the aforementioned press, when the sheet dries are reproduced [on paper] as more contracted in one direction. 
But neither our Riccioli\footnote{Giovanni Battista Riccioli SJ (1598--1671), Italian astronomer, author
of \emph{Geographiae et hydrographiae reformatae libri duodecim}, Bologna, 1661.} in his ``Geography reformed'' avoided all dangers of error by arranging printing of the Roman half-foot
on dried sheet\footnote{See Fig. \ref{semipesfig} in the Appendix.}; for even this sheet to some extent changes with changes of atmospheric conditions; this is not much different than how the width of wooden tables or even ivory diminishes with time.
A more secure method, therefore, of transmitting particular measures to posterity
would be to engrave them in the metal form, and then impress on gold leafs, of the kind used
by garland makers, called \emph{Orpello} in Italian or \emph{Rauchgold} in German. These could then
be easily inserted into books by a pressman. The moderate changes which they undergo subjected to heat or cold do not induce any detectable error.
If somebody wants to obtain genuine measures from such imprints alone, one can reproduce them
accurately enough  from a sheet previously moistened and suitably extended.

Measure of half or quarter of foot can be most conveniently mailed to a friend by letter  
as a piece of bronze wire, of the kind used in musical instruments.
On its ends two loops are made by twisting, and the twisted parts are crushed a bit with the hammer
so that they better hold together. Finally, stretching  the bronze wire moderately
over the given measure, its internal ends are accurately trimmed with a small knife or scissors, as much as needed. Inspect the provided figure.
\pend

\end{Rightside}

\Columns

\end{pairs}
\begin{center}
\includegraphics[width=16cm]{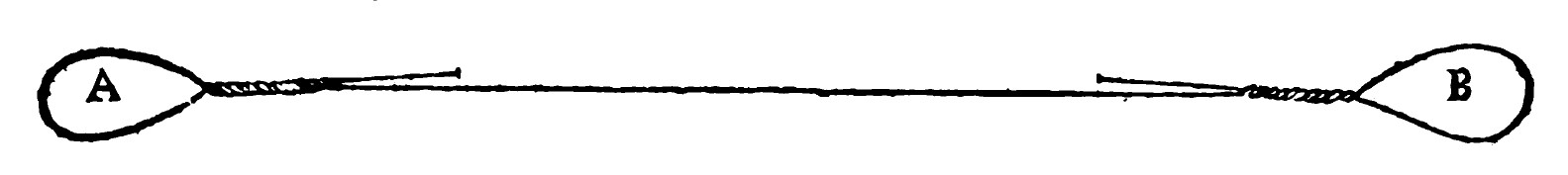}
\end{center}
\begin{pairs}
\begin{Leftside}
\selectlanguage{latin}
\noindent\setline{1}
\pstart
IV. \emph{Penduli portatilis, ac Horologiorum perfectio.}

Vibrationes ejusdem penduli in arcubus in aequalibus inaequales
esse, a multo tempore experientia oculis ipsis manifesta 
didici. Hanc quidem inaequalitatem Cycloidis adminiculo 
sublatam esse non ignoro, videtur interim mihi Cycloidis usus cum
pendulo e fascia appenso, in usus nauticos, ad longitudines 
Terrae investigandas, ob navis agitationem, nonnihil impeditus 
esse. Quamobrem construxi horologium, cujus rotae ad numeros
Astronomicos decurrunt, in eoque duo pendula aequalia 
appensa sunt ex elateribus chalybeis, mediocriter latis, ac facile 
flexibilibus, \&  breviusculis quidem, sed qui pro re nata produci
nonnihil, vel contrahi possunt, reliquae pendulorum virgae sunt
rigidae. Experior autem elateres hos adinstar cycloidis, digressiones
\&  in iis moras pendulis non permittere, sed aequales ad 
sensum vibrationes peragere: quae quidem in horologio meo 
breves sunt. Nam binae compositae five quaternae simplicas uni 
secundo respondent: verum \& in longioribus pendulis e 
convenienti elatere  pendentibus, de felici successu dubitare non 
licet. \emph{Experiantur itaque navigantes}, annon pendulum elatere
hujusmodi sustentatum, altero illo suspenso e fascia cycloidem
stringente, sit in mari commodius.
\pend
\pstart
Horologiis rotatis in sacco portatilibus applicueram olim 
pendulum Magneticum, quod frustro validi magnetis ad vibrationes
incitabatur satis aequales; sed quia pondus lapidis gravat, \&
ferri vicinia aequalitatem vibrationum perturbat, ideo repudiato
hoc invento, id tandem effeci, ut horologia parva nunc usitata,
quae librile habent elatere spirali instructum, nulla inter
ambulandum vel equitandum agitatione turbentur \emph{in suo} motu
regulari. Appendo enim totum horologii corpus intra \emph{suam}
thecam extimam nonnihil ampliorem, ita ut circa systematis
totius axem, librilis  axi parallelum, vel circa duos potius eius
polos, aut cuspides, instar Rosae nauticae, libere convertatur
horologium totum: quo fit, ut nulla agitatio horologii sic 
appensi ejus librile in suo cursu perturbet, quemadmodum sit in
usitatis , quae manu velociter reciproceque agitata, vibrationes
librilis enormiter praecipitant: Idipsum autem etiam in accelerata
deambulatione, quando videlicet sacciperium pendens, una
cum horologio in eo contento, reciproce agitatur, contingere
solet, \& incautos decipere. Praxis appensionis illius varie
poterit institui. Ego hanc tenui. Thecam extimam ansa \& annulo
appensorio, nec non crystallo instruxi: in fundo thecae
defixi brevem stylum, qui pileolum in pixide horologii firmatum
sustentat: alter polus est ille ipse, qui indicem horologum 
defert, \&  in altero pileolo crystalli centrum occupante,
circumagitur.
\pend
\pstart
 Non ingratum fortasse Curiosis erit hac occasione intelligere,
me excogitasse ab annis aliquot artificium quoddam, quo
mediante, quilibet rudis homo, puer aut foemina, quin \&  
caecus incidere poterit quascunque minutas horologiorum rotas,
tot dentium quot imperantur, absque praevia earum divisione,
\& insuper eadem opera plures, aequales, maiores, \&  minores.
Praeterea aliud artificium ad formandam Torno cochleam 
illam Horologii, quae chordam, vel catenulam colligit, dando
illi figuram, attemperandis viribus cujuscunque dati elateris
chalybei convenientem; ita ut cochlea circuitus imperatos 
complectatur. An autem ista publicari debeant, ne in eorum artificum,
qui veteribus insistere coguntur, detrimentum cedant,
anceps haereo.

  Ausim etiam Opticis polliceri Tornum universale, in quo
lentes , \&  specula omnium sectionum Conicarum, elaborari
possunt; \&  quidem ratione modoque inter omnes alios 
commodissimo: reddidi enim illud immune eo defectu, quo omnes
aliorum machinae tales laborant; ne videlicet inter laborandum,
propositae figurae modulus deteratur.

  Hucusque ea, quae a me in Galliam missa fuerant, e quibus
illud posterius Num. IV de Horologiis in sacco portatilibus,
quia jam a me in Actis Erudit. Anni 1685 Mense Septembri
commemoratum fuerat, ne hoc loco cramben bis coctam 
sapiat, hoc ei condimentum adjiciam. Posse videlicet, absque
nova totius horologii constructione, qualem praedicto loco 
supposui , omnibus aliis Portatilibus, illam circa binos axes 
conversionem, non difficulter conciliari, ac praesertim iis, quae
spirali, aut etiam recto elatere vibrationes librilis uniformes
efficere debebunt. Horologium itaque usitata ratione jam 
constructum, thecae cuidam extimae (rejecta vel retenta, si lubet,
pristina) includatur, intra quam aptata fit armilla quaedam 
latior, in Ellipsin, horologio recipiendo, \&  firmiter retinendo
idoneam, inflexa, \&  in sua minore Diametro, duabus oppositis
cuspidibus instructa, circa quas, tanquam Polos, totum 
Horologium cum armilla intra thecam extimam, armilla nonnihil
ampliorem, converti circulariter possit: In Horologiis enim
taliter gyrantibus, eorum librili concedi poterit elater fortior,
ac praedominans elateri majori, rotas omnes impellenti, ut hoc
modo vibrationes librilis, inaequali rotarum dentatarum tractioni
minus obnoxia sint, nec tamen primum impulsum sibi datum,
\& a totis incitantibus conservatum, amittant.
\pend
 
\end{Leftside}

\begin{Rightside}
\selectlanguage{english}
\noindent\setline{1}
\pstart
IV. \emph{Perfection of portable pendula and watches.}
I have learned, with experiments revealing it  to the eyes many times, that oscillations of
the same pendulum in inequal arcs are inequal.\footnote{In modern terminology:
a pendulum's period of oscillations depends on the amplitude.}
I am not ignorant that cycloidal support removes this inequality\footnote{Cycloidal pendulum was already known to be isochronous, as 
demonstrated by Christiaan Huygens in 1673.}, yet it seems to me that
the use of Cycloid with the pendulum suspended on a ribbon in the nautical service,
for the purpose of determining the longitude, would be somewhat impedited by the 
movement of ship. Because of this I constructed a clock whose wheels run
according to astronomical numbers, in which two equal pendula hang from steel springs,
moderately wide and easily flexible yet very short, which by nature can stretch or
contract a little bit, with the remaining rods of pendulums being rigid.
I attest that these springs, similarly to cycloids, do not permit the pendula to
run into delays but force them into perceptibly equal  oscillations. 
Two composite or four simple [oscillations] correspond to one second: but certainly
even with longer pendula, suspended from suitable springs, one should not doubt
of successful outcome. Let, therefore, navigators test which one is more desirable
on the sea, a pendulum supported by the spring this [described] way, or
another one suspended on a ribbon restricted by the cycloid.
\pend
\pstart
I once adopted to portable watches a magnetic pendulum, induced to equal enough vibrations
by a strong piece of magnet; but because the weight of the stone, and because the vicinity of iron
perturbs the oscillations, abandoning this invention, I finally made it so that 
the small watches now used, which have the balance equipped with a spiral spring, 
are not perturbed in their regular motion by walking or horse riding.
I suspend the entire body of the watch inside of an external slightly larger case,
so that the entire watch rotates freely around the axis of the entire system, parallel to
the axis of the balance, or rather around two its poles or pointed ends, similar to the nautical
rose. For this reason it happens that no shaking of the watch suspended this way perturbs
the balance from its course, in whatever way the watch is used, and if it is  violently
moved back and forth with the hand, the oscillations of the balance are hugely cast down.
Moreover, this happens  and surprises observers even in an accelerated walk, when the hanging sack evidently shakes
back and forth together with the  watch kept in it.
In practice this suspension can be set up in various ways. Here is how I do it.
I constructed the outer case with a hook and loop handle, as well as crystal.
In the bottom of the case I attached small pillar which supports a cup in the 
small case of the watch. Another pole is the one which carries the hand of the watch,
and is rotates held in place by another cup in the center of the crystal.
\pend
\pstart
Those eager to know will perhaps not be ungrateful to learn on this occasion that 
some years ago I devised a certain method, by which whichever simple man, boy or woman,
in fact even a blind one will be able to cut however small clock wheels, with as many teeth
as ordered, without dividing them beforehand, and by the same effort make many of them,
equal ones, smaller ones and bigger.
In addition to this, [I devised] another method for forming on the lathe 
the clock's snail, which holds a rope or a little chain, giving it its appropriate shape,
with adjustable strength of a given steel string.
I am still undecided whether this should be published, so it does not become lost, like those skills
which are restricted to be practiced by the old only.

I would also dare to promise to opticians  a universal lathe on which one can make all lenses and mirrors with conical surfaces; its plan and the way of operation are most agreeable among all others.
For I rendered it immune to this defect with which all other machines struggle: namely,
while operating, it does not wear away the dimension of the proposed shape.

In order for the observations about portable watches in a bag,
which were sent by me to France, mentioned in section IV  and placed on record in
Acta Euruditorum in September 1685, not to taste like doubly-cooked cabbage, 
I will add some condiments here.
Namely, it is possible, without difficulty, to procure the conversion of all other portable
watches, without construcing the whole new watch with the aforementioned method,  so that they rotate around two axes:  this especially applies to
watches with balances which need to oscillate uniformly using a spiral or even a straight  
spring.
Let  the watch ordinarily constructed be included in a case (original case
discarded or retained, as preferred), inside which  
a ring is fitted, bent in the shape of ellipse, accepting the watch, and holding it firmly
in a suitable way; let its minor diameter be equipped with two opposite 
pointed ends. Around those ends, as if they were poles, entire watch with the ring
can rotate inside the outer case.
In watches gyrating this way, a stronger spring could be subjected to their
balance, and dominate over the major spring, driving all wheels. This way  
oscillations of the balance would be less liable to the uneven pull of gears,
and less liable to loosing original momentum given to it, which is also protected from 
sudden increases.
\pend

\end{Rightside}

\Columns

\end{pairs}

\newpage
\section*{Appendix}
\begin{figure}[h!]
 \begin{center}
  \includegraphics[width=10cm]{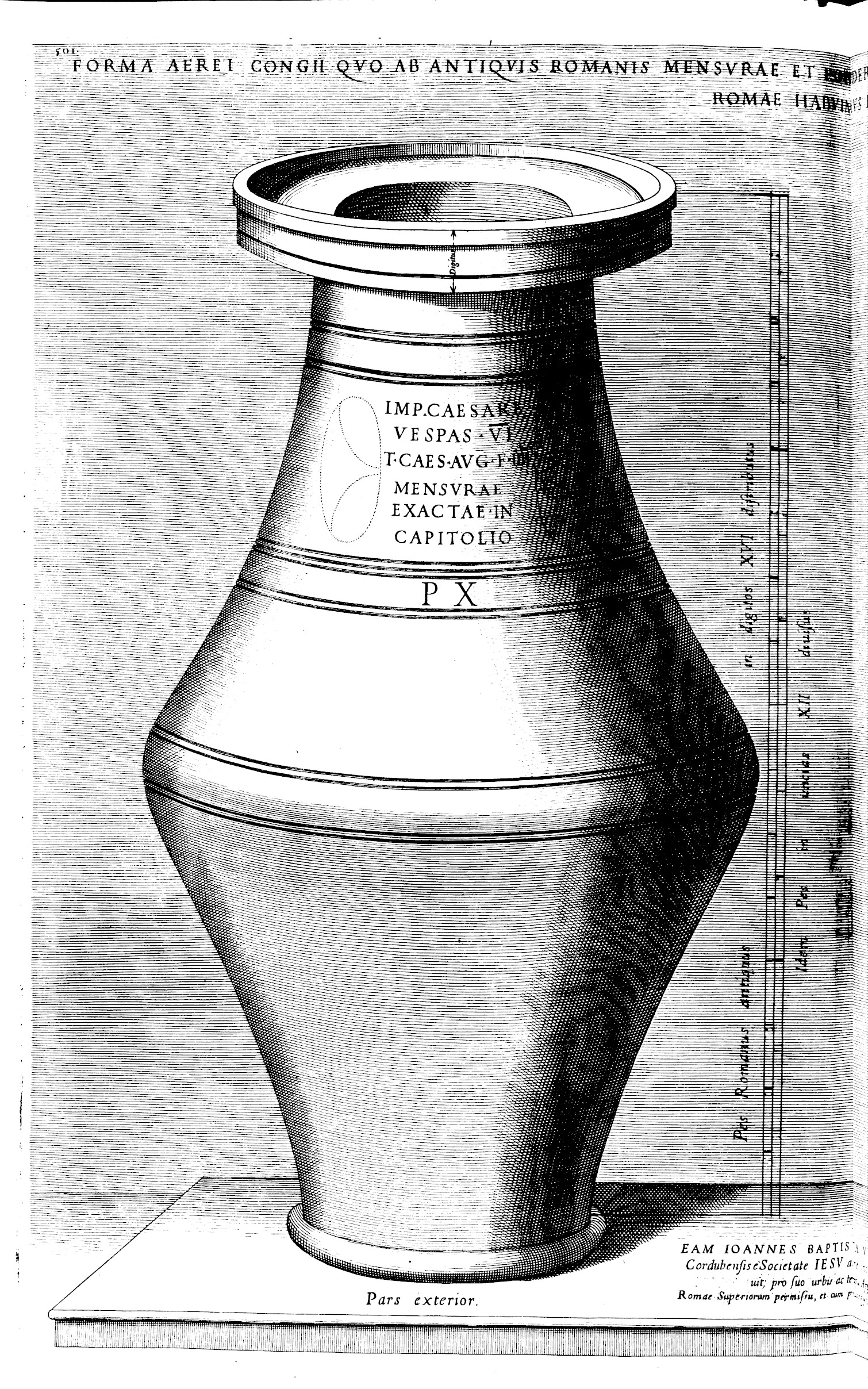}
  \caption{Drawing of the Roman \emph{congius} togther with
  the length of the Roman foot, reproduced from J. B. Villapando, \emph{Apparatus Urbis ac Templi Hierosomitani}, Romae 1604.}
  \label{vasefig}
 \end{center}
\end{figure}
\begin{figure}[h!]
 \begin{center}
  \includegraphics[width=13cm]{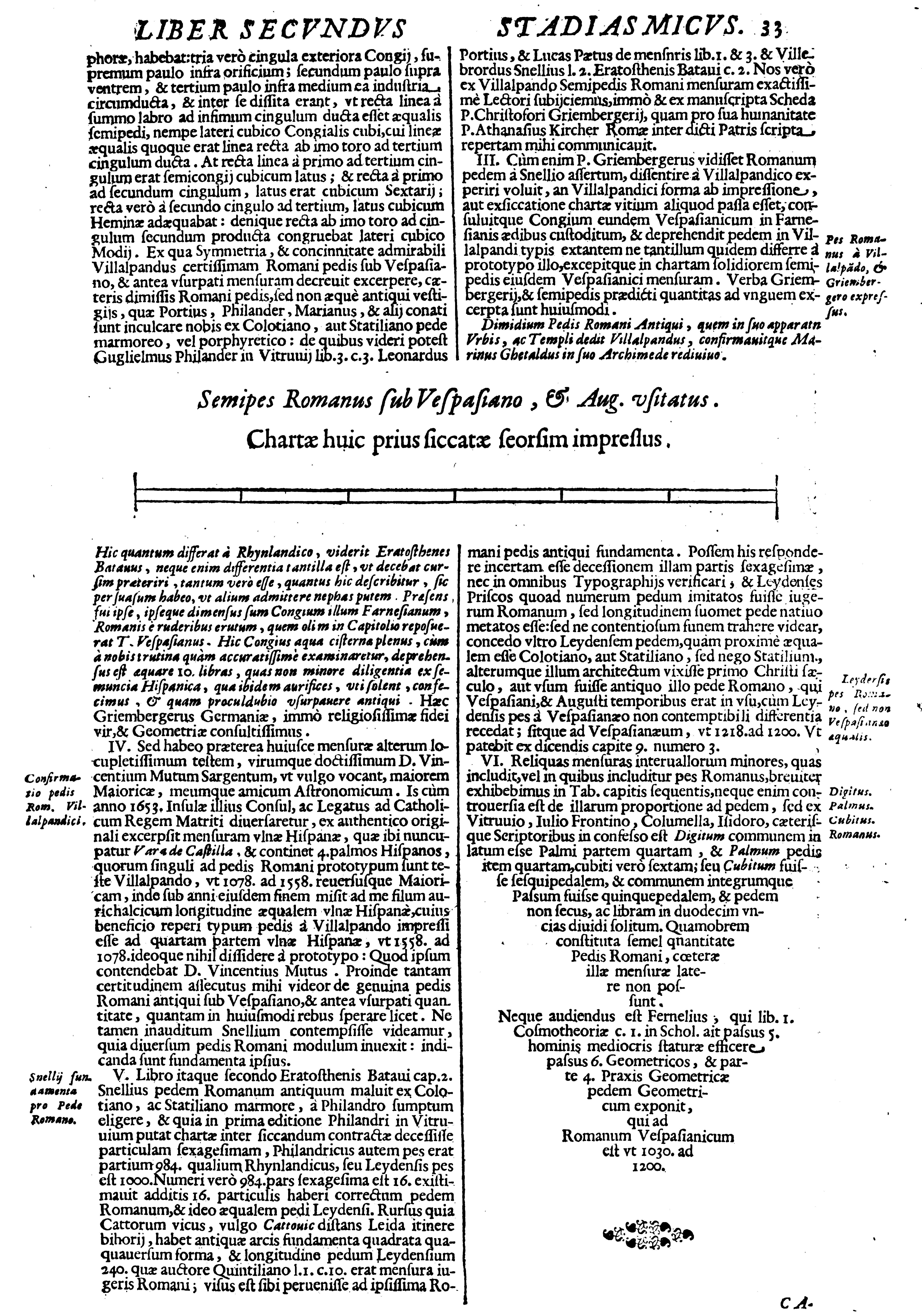}
  \caption{Length of the Roman half-foot, from
  G. B. Riccioli, \emph{Geographiae  et hydrographiae reformatae libri duodecim}, Bologna, 1661.  The caption above says
  ``Separately printed on this sheet previously  dried''.}
  \label{semipesfig}
 \end{center}
\end{figure}

\end{document}